\documentclass[paper=a4, fontsize=12pt]{scrartcl} % , twocolumn KOMA-article class

\usepackage{comment}
\usepackage[english]{babel} % English language/hyphenation
\usepackage[protrusion=true,expansion=true]{microtype} % Better typography
\usepackage{amsmath,amsfonts,amsthm} % Math packages
\usepackage{enumitem} %allows to change enumeration style 
\usepackage[pdftex]{graphicx} % Enable pdflatex
\usepackage[svgnames]{xcolor} % Enabling colors by their 'svgnames'
\usepackage[hang, small,labelfont=bf,up,textfont=it,up]{caption}	% Custom captions under/above floats
\usepackage{subfig}	% Subfigures
\usepackage{fix-cm} % Custom fontsizes

\usepackage{url}
\usepackage{marginnote}

\usepackage{sectsty} % Custom sectioning (see below)
\allsectionsfont{%															% Change font of al section commands
	\usefont{OT1}{phv}{b}{n}%										% bch-b-n: CharterBT-Bold font
}

\sectionfont{%																% Change font of \section command
	\usefont{OT1}{phv}{b}{n}%										% bch-b-n: CharterBT-Bold font
}

\usepackage{fancyhdr} % Needed to define custom headers/footers
\usepackage{lastpage}	

\fancypagestyle{default}{
	\lhead{}
	\chead{}
	\rhead{}
	\lfoot{\footnotesize Petra Schwer \textbullet ~How to Read Mathematics? \textbullet ~Licence CCBY1 \textbullet ~2026}
	\cfoot{}
	\rfoot{\footnotesize page \thepage\ of \pageref{LastPage}}	% "Page 1 of 2"

}

\fancypagestyle{last-page}{
	\lhead{}
	\chead{}
	\rhead{}
	\lfoot{\footnotesize Petra Schwer, Professor of Mathematics at Heidelberg University, Germany. \\
		Share your thoughts with me at  \emph{schwer@uni-heidelberg.de} \\
		\url{https://web.mathi.uni-heidelberg.de/ggt/schwer} \\ ~ \\
		How to Read Mathematics? \textbullet ~2026 \\
		This text is licenced under Creative Commons Licence CCBY1. }
	\cfoot{}
	\rfoot{~ \\~ \\~ \\~ \\
		\footnotesize \includegraphics[width=15ex]{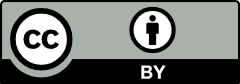}	 } %page \thepage\ of \pageref{LastPage}}	% "Page 1 of 2"

}

\AtEndDocument{\thispagestyle{last-page}}% Page style at \end{document}

\usepackage{lettrine}
\newcommand{\initial}[1]{%
\lettrine[lines=3,lhang=0.3,nindent=0em]{
\color{DarkGoldenrod}
{\textsf{#1}}}{}}

\usepackage{titling}% For custom titles

\newcommand{\HorRule}{\color{DarkGoldenrod}	\rule{\linewidth}{1pt}}

\pretitle{\vspace{-30pt} \begin{flushleft} \HorRule 
\fontsize{40}{40} \usefont{OT1}{phv}{b}{n} \color{Maroon}%{DarkSlateGray}%{DarkRed} 
\selectfont 
}

\title{\Huge{How to Read Mathematics? }\\
\Large{{\ }\\ A study guide.
	}
} % Title of your article goes here
\posttitle{\par\end{flushleft}\vskip 0.5em}

\preauthor{\begin{flushleft} 
	\large \lineskip 0.5em \usefont{OT1}{phv}{b}{sl} \color{black}}%{DarkRed}}
	\author{Petra Schwer, }											% Author name goes here
	\postauthor{\footnotesize \usefont{OT1}{phv}{m}{sl} \color{Black} 
Heidelberg, Summer 2026
\par\end{flushleft}\HorRule}
\date{}																				% No date

\begin{document}
\maketitle
\thispagestyle{fancy} 			% Enabling the custom headers/footers for the first page 

% The first character should be within \initial{}
\initial{R}\textbf{eading a math textbook isn't like reading a novel. You've probably already figured that out on your own. But how *do* you read mathematical texts? And how can you actually learn from a book or lecture notes?} 

This document is a guide to independently reading and working through mathematical texts.  
With the methods and steps explained here, you can work through any paper, textbook, lecture notes, or even your own class notes.  
Reading math is hard. But don’t worry -- hard things are exactly where growth happens!

%%%%%%%%%%%%%%%%%%%%%%%%%%%%%%%%%%%%%%%%%%%%%%%%%
\subsection*{\textcolor{Maroon}{General Tips}}

Let me start with a few foundational tips for reading mathematical texts. Just mastering these six points, once internalized and used regularly, can make a huge difference.

\bigskip

\noindent
\textbf{Read slowly.}

No, even slower. One does not read math in pages per hour but in hours per page.

\medskip
\noindent
\textbf{Give yourself time.}

Be prepared: it can take a long time to understand a single page, and in most cases, you’ll need more than one try. Plan to work on it over several days. At the same time, don’t make your study sessions too short. Math needs repetition and time. Time to work with the material, and time for the ideas to simmer in your mind.

\medskip
\noindent
\textbf{Ask precise questions to yourself and others.}

This one’s crucial: questions are the key to deep understanding.  
"I don’t understand Definition 3.12" isn’t a question.  
"How should I interpret Definition 3.12?" isn’t precise enough.  
"Why must the inverse function itself be bijective in the definition?" now *that’s* a good question!

\medskip
\noindent
\textbf{Talk about the material with others.}

Try to express yourself as clearly as possible and use proper mathematical language. Even if you think you’ve got it, talk through it anyway. Sometimes you only realize you don’t truly understand something when you try to explain it to someone else.

If you’re stuck for too long, don’t hesitate to ask others or consult different sources. But don’t just leave them sitting side by side. Integrate those sources into a coherent picture. One great way to do this? Write it down using uniform notions and notation. Which leads me to the following.  

\medskip
\noindent
\textbf{Reading is also writing.}

Highlight confusing parts. Write summaries for each topic or key proofs (not every lemma). Write down every detail of the examples. The latter are sometimes kept brief in textbooks. Write out all the reasoning steps. Come up with your own examples and write them down too.

\bigskip

The next two tips might seem contradictory but both are important. Find your own balance between:

\medskip
\noindent
\textbf{Keep going vs. Move forward}

Keep going: Don’t give up the moment you don’t understand a line. Try to figure it out.  
Move forward: If something truly isn’t clicking, keep reading and come back to it later. Often, a clearer explanation a few lines down will help for things to make sense.

%%%%%%%%%%%%%%%%%%%%%%%%%%%%%%%%%%%%%%%%%%%%%%%%%
\subsection*{\textcolor{Maroon}{How to Read a Mathematical Text?}}

In three (or four, or five) steps.

\medskip
\noindent
Step one can happen anywhere: on the couch, on the bus, on a park bench in the sun… wherever you can hold your book or device. My son loves reading on the swing. Whether or not you take notes in this pass is up to you.

\subsubsection*{1. Get a quick overview}

Skim over the text. You can read it casually, like a regular book. The goal is just to get a sense of what’s in there. Ask yourself:

\begin{enumerate}[label=\alph*]
\item What is this text about? How is it divided into logical sections? These logical sections (which may or may not agree with actual sections) can become your work packages for later. No need to understand everything at once.
\item What are the main ideas?
\item Which part looks tricky or especially hard and why?
\item Are there lots of new definitions and concepts? Or many long proofs?
\item What do I already know about this topic? Does anything in this chapter feel familiar?
\end{enumerate}

Now you’ve got a rough sense of the material. You’ve probably already started memorizing key terms, names of objects, methods, or theorems. Time to dig deeper and actually understand the content. Ready?

\medskip

For steps two and three, find a quiet place where you can focus without distractions. Thirty minutes of uninterrupted work beats an hour of distracted reading.

\subsubsection*{2. Understand definitions and statements}

You are now going to build deeper understanding. Remember: this takes time. Here are some questions to ask yourself while reading.

\bigskip

\noindent
\textbf{Reading definitions:}

Definitions are the vocabulary of math. Without precise use, we couldn’t communicate clearly about math at all.  
Examples exist to give you intuition about how and where a definition applies and where it doesn’t.

For every definition you read, ask:
\begin{itemize}
\item What is being defined here? Is it a new property of something familiar, or a completely new object?
\item Why is this concept well-defined?
\item How important is this idea? Is it just a technical tool for a proof, or a central concept in the theory? (This isn’t always obvious at first.)
\end{itemize}

No definition without examples:
\begin{itemize}
\item What’s an example of this definition?
\item What’s a non-example? Can I find an object that looks similar but fails one or more conditions?
\item Can I follow the given examples? Can I fill in any missing logical gaps myself? If not, go back and reread the definition.
\end{itemize}

\medskip

\noindent
\textbf{Reading lemmas/theorems/corollaries:}

Lemmas and theorems state mathematical truths about how objects defined earlier behave. They’re like character sketches or short stories about how things relate.

For every mathematical statement, ask:
\begin{itemize}
\item Do I understand the statement? That means: Can I identify the assumptions and the conclusion? Can I make sense of the conclusions?
\item What are some examples of this statement? How does it look in a concrete case?
\item What happens if one of the assumptions fails? Does the statement still hold? If not, what breaks?
\item Why is this a lemma/theorem/corollary?
\item What does this statement explain, and how does it fit into the bigger picture?
\end{itemize}

Also read remarks (which offer extra explanations, intuition, or ideas), introductory or transitional text, and notation sections (where writing and wording conventions are set). Don’t skip these as they’re part of the story.

Step three usually takes the longest. But you can break it into small chunks. If you’re stuck on a proof, it’s okay to come back later after working through step 4.

\subsubsection*{3. Understand the proofs}

Now it gets exciting. Go deep and figure out *why* something works the way it does. Let’s read proofs!

Proofs aren’t just about verifying a statement. They help you understand connections and properties of mathematical objects.

Treat each proof like a mini-mathematical text. Start again with step one: get an overview!

\begin{enumerate}[label=\alph*]
\item Skim the proof.
\begin{itemize}
	\item What proof technique is being used?
	\item Are there intermediate claims or sub-statements?
	\item What’s the main idea? (This isn’t always obvious at first. If you don’t see it yet, keep reading.)
	\item Which parts look difficult?
	\item Where are the assumptions used? (This can be tricky. Keep reading if it’s unclear.)
\end{itemize}

\item Read the proof carefully, line by line.
\begin{itemize}
	\item Do I understand each step and transformation?
	\item What assumptions or properties of the objects do I need at each step?
	\item Have I seen a similar proof before?
	\item Can I walk through the proof using a concrete example?
\end{itemize}
\item Summarize the key points.
\begin{itemize}
	\item What are the main ideas of the proof?
	\item Where are the assumptions used?
\end{itemize} 	
\end{enumerate}

\bigskip

If you’ve made it this far, you’ve already done a lot of hard work and accomplished a lot. Now you’re probably wondering: How do I know I actually *understand* what I’ve read? How can I check my understanding?

The truth is: you learn math by doing. You don’t have to invent new math but you do need to make the new math your own. The most powerful tool for this? Step four:

\subsubsection*{4. Work through practice problems}

You could write a whole guide just on this topic. But here are the essential steps:

\begin{enumerate}[label=\alph*]
\item Look at the problem set as soon as you get them. \newline That way, the problems are already in your head, and you can start thinking about them. Treat the problem sheet like another math text and use the tips from steps 1–3 to read them.
\item Understand the problem statement. \newline Yes, this is a separate step. You can’t start solving unless you know what you’re being asked to do. If you’re stuck, ask: Do I still remember and understand the relevant definitions?
\item Talk about the problems with others. \newline This is done best with people working on the same problems. Brainstorming together helps.
\item Work on the problem with focus. \newline Find a quiet space. Think carefully: Which definitions from class are relevant? Do I know them? Which theorems seem related and might help? Can I imitate a proof we've seen in class?
\item Write down your solution. \newline Just having the idea in your head isn’t enough. Writing math clearly is a separate skill which is sometimes harder than coming up with the idea. Even if you’re sure it’s perfect, write it down.
\item Explain your solution to someone else. \newline You can do this even if you’re not 100\% confident. Ask others to spot errors or double-check your logic. No one’s going to eat you alive.  
\end{enumerate}

\subsubsection*{5. Rebuild the big picture}

At the very end of the reading process, go back and repeat step one. Rebuild your understanding in a broader context. This helps you avoid getting lost in the details of definitions, examples, lemmas, theorems, and proofs. You know the saying: “You can’t see the forest for the trees.”

Now, take a structured look at the whole topic. This time, take notes. Ask yourself the same questions from the beginning:

\begin{enumerate}[label=\alph*]
\item What is this text about?
\item What are the main ideas?
\item What are the most important definitions and statements?
\item What still feels hard? What do I feel confident about?
\end{enumerate}

Sometimes, even after all of this, you’ll still have unanswered questions. But many things will already be clearer. Go back to step one. Talk to others. Be persistent. Things will start to click. Sometimes right away, sometimes much later.
Either way, be proud of yourself. You’ve worked hard and learned a lot.\newline Congratulations!

%%%%%%%%%%%%%%%%%%%%%%%%%%%%%%%%%%%%%%%%%%%%%%%%%
\subsection*{\textcolor{Maroon}{The Inevitable Note on AI}}

The role of AI in education is a huge topic. Too big to cover it in full here.

\bigskip

Students and researchers can gain a lot from the tools and insights AI provides. And one needs to use the results with caution and always verify correctness.

The danger lies in relying too much on AI and losing your own ability to think. Math is challenging. Understanding new concepts can take time. Getting comfortable with definitions, theorems, and their use can be tough. But that’s exactly how your brain builds new connections.

These observations aren’t surprising or new and definitely not just mine.

\bigskip

The reading methods outlined above don’t mention AI. In my view, though AI might help sometimes, it’s not necessary and could even be harmful in some cases.
Decide for yourself when and how to use AI. Consider the impact of using it. Here are two related findings from recent studies:

\bigskip

%\emph{Using EEG data, researchers tracked how each group engaged cognitively [...]. The results were striking. Students using ChatGPT showed lower brain activity, weaker memory recall, and less ownership of their writing.} \footnote{From "Thinking with AI: What a viral MIT study reveals about AI and the learning brain" by Tina Peters, Risepoint. % \url{https://risepoint.com/insights/thinking-with-ai-what-a-viral-mit-study-reveals-about-ai-and-the-learning-brain/};}

\emph{EEG revealed significant differences in brain connectivity: Brain-only participants exhibited the strongest, most distributed networks; Search Engine users showed moderate engagement; and LLM users displayed the weakest connectivity. [...] %Cognitive activity scaled down in relation to external tool use. In session 4, LLM-to-Brain participants showed reduced alpha and beta connectivity, indicating under-engagement. Brain-to-LLM users exhibited higher memory recall and activation of occipito-parietal and prefrontal areas, similar to Search Engine users. Self-reported ownership of essays was the lowest in the LLM group and the highest in the Brain-only group. LLM users also struggled to accurately quote their own work. 
	While LLMs offer immediate convenience, our findings highlight potential cognitive costs. Over four months, LLM users consistently underperformed at neural, linguistic, and behavioral levels. These results raise concerns about the long-term educational implications of LLM reliance and underscore the need for deeper inquiry into AI's role in learning. }\footnote{Aus: "Your brain on chatgpt: Accumulation of cognitive debt when using an ai assistant for essay writing task." von Nataliya Kosmyna, Eugene Hauptmann, Ye Tong Yuan, Jessica Situ, Xian-Hao Liao, Ashly Vivian Beresnitzky, Iris Braunstein, and Pattie Maes. preprint 2025, \url{https://arxiv.org/abs/2506.08872}}

\bigskip

\emph{The findings reveal that AI offers significant benefits, including personalized learning, improved academic outcomes, and enhanced student engagement. However, challenges such as over-reliance on AI, diminished critical thinking skills, data privacy risks, and academic dishonesty were also identified.} %The study underscores the necessity of a structured framework for AI integration, supported by ethical guidelines, to maximize benefits while mitigating risks.}
\footnote{Aus: "The Impact of Artificial Intelligence (AI) on Students’ Academic Development",  von Aniella Mihaela Vieriu and Gabriel Petrea, Educ. Sci. 2025, 15(3), 343; \url{https://doi.org/10.3390/educsci15030343}}

\end{document}